\documentclass[11pt]{amsart}

\usepackage[T1]{fontenc}
\usepackage[utf8]{inputenc}
\usepackage{amsmath,amssymb,amsthm,mathtools}
\usepackage{tikz-cd}
\usepackage{xcolor}
\usepackage{microtype}
\usepackage{booktabs}
\usepackage{url}
\usepackage[margin=3cm]{geometry}
\usepackage[pagebackref]{hyperref}
\hypersetup{
  colorlinks=true,
  linkcolor = violet,
  citecolor = blue,
}
\renewcommand*{\backref}[1]{}
\renewcommand*{\backrefalt}[4]{%
  \ifcase #1\relax
  \else\space$\uparrow$~#2%
  \fi}
\usepackage[nameinlink,noabbrev,capitalise]{cleveref}
\newcommand{\C}{\mathbb{C}}
\newcommand{\Z}{\mathbb{Z}}
\newcommand{\Q}{\mathbb{Q}}

\newcommand{\PP}{\mathbb{P}}
\newcommand{\OK}{\mathcal{O}}

\DeclareMathOperator{\vol}{vol}
\newcommand{\volhat}{\widehat{\operatorname{vol}}}

\DeclareMathOperator{\lct}{lct}
\DeclareMathOperator{\codim}{codim}

\DeclareMathOperator{\ord}{ord}

\DeclareMathOperator{\Sym}{Sym}
\DeclareMathOperator{\Supp}{Supp}

\DeclareMathOperator{\rk}{rk}
\DeclareMathOperator{\im}{im}

\DeclareMathOperator{\PGL}{PGL}

\theoremstyle{plain}
\newtheorem{theorem}{Theorem}[section]
\newtheorem{proposition}[theorem]{Proposition}
\newtheorem{lemma}[theorem]{Lemma}
\newtheorem{corollary}[theorem]{Corollary}
\theoremstyle{definition}
\newtheorem{definition}[theorem]{Definition}
\newtheorem{conjecture}[theorem]{Conjecture}
\theoremstyle{remark}
\newtheorem{remark}[theorem]{Remark}
\crefname{section}{Section}{Sections}
\crefname{subsection}{Section}{Sections}
\Crefname{section}{Section}{Sections}
\Crefname{subsection}{Section}{Sections}

\title[K-stability and GIT for cubic hypersurfaces]
      {K-stability and GIT for cubic hypersurfaces}

\author{Long Pan}
\address{Academy of Mathematics and Systems Science, Beijing, China}
\email{panlong@amss.ac.cn}

\author{Haoyu Wu}
\address{Yau Mathematical Sciences Center, Tsinghua University, Beijing, China}
\email{haoyu-wu@tsinghua.edu.cn}

\date{\today}



\begin{document}

\begin{abstract}
	In this note, we prove that K-(semi/poly)stability of cubic hypersurfaces of dimension at least five coincide with the corresponding GIT stability.
	In particular, their K-moduli and GIT moduli stacks are isomorphic.
	The main step is to show that every K-semistable $\Q$-Gorenstein limit of smooth cubic hypersurfaces is again a cubic hypersurface.
	We establish Cartierness of the limiting hyperplane divisor using local log canonical thresholds and mixed intersections.
\end{abstract}

\maketitle

\begingroup
\small\linespread{0.90}\selectfont
\tableofcontents
\endgroup

\section{Introduction}
\label{sec:introduction}

\subsection{Background and main result}
\label{subsec:problem}

K-stability, introduced by Tian \cite{Tian97} and reformulated algebraically by Donaldson \cite{Dona}, relates the existence of K\"ahler--Einstein metrics to the geometry of degenerations of Fano varieties.
The algebraic theory now provides moduli stacks of K-semistable Fano varieties with projective good moduli spaces parametrizing K-polystable varieties; see \cite{Reductivity,Openness,Blum-Xu-uniqueness,CM,Jiang,LWXb,Liu-Xu-Zhuang-Finitegeneration} and the references therein.
Describing these moduli spaces explicitly, however, requires control of the varieties that occur on their boundary.

For several families, the moduli continuity method identifies the K-moduli compactification with a classical GIT quotient.
Examples include del Pezzo surfaces \cite{OSS}, intersections of two quadrics \cite{ODP}, cubic threefolds \cite{Liu-Xu-cubic-threefolds}, and cubic fourfolds \cite{Liu-cubic-fourfolds}.
A central ingredient is the comparison between the anticanonical volume of a K-semistable Fano variety and the normalized volumes of its singularities.
The latter invariant was introduced by C.~Li \cite{LiChi-normalizedvolumes}.
For varieties of large anticanonical volume, this comparison places strong restrictions on possible degenerations.

In this paper we study cubic hypersurfaces in higher dimensions.
If $X\subset\PP^{n+1}$ is a smooth cubic $n$-fold and $H=\OK_X(1)$, then
\[
	-K_X\sim(n-1)H,\qquad H^n=3.
\]
Liu--Xu \cite[Theorem~1.1]{Liu-Xu-cubic-threefolds} and Liu \cite[Theorem~1.1]{Liu-cubic-fourfolds} proved that K-stability and GIT stability agree in dimensions three and four, respectively.
Kong--Shen--Zhao--Zhao \cite{Kong-Shen-Zhao-Zhao-cubic-fivefolds} established the corresponding result for cubic fivefolds.
We prove the following statement uniformly for all $n\geq5$.

Let $\mathcal{M}^K_{\mathrm{cubic}}$ denote the K-moduli stack parametrizing K-semistable Fano $n$-folds of anticanonical volume $3(n-1)^n$ that admit a $\Q$-Gorenstein smoothing to a smooth cubic.
Let $U^{\mathrm{ss}}$ be the GIT-semistable locus in the projective space of cubic forms on $\PP^{n+1}$, and set
\[
	\mathcal{M}^{\mathrm{GIT}}_{\mathrm{cubic}}
	= [U^{\mathrm{ss}}/\operatorname{PGL}(n+2)].
\]
We test GIT stability using $\PGL(n+2)$ and use the effective projective action for the quotient stack.

\begin{theorem}\label{thm:main}
	Let $n\geq5$, and let $X=V(f)\subset\PP^{n+1}$ be a cubic hypersurface.
	Then $X$ is K-semistable (respectively, K-polystable or K-stable) if and only if the cubic form $[f]$ is GIT-semistable (respectively, GIT-polystable or GIT-stable) for the natural action of $\PGL(n+2)$.
	Moreover, there is a natural isomorphism
	$\mathcal{M}^K_{\mathrm{cubic}}\cong
		\mathcal{M}^{\mathrm{GIT}}_{\mathrm{cubic}}$.
\end{theorem}

\begin{proof}
	The CM polarization of the universal cubic family is a positive rational multiple of the GIT polarization.
	Thus K-(semi/poly)stability implies the corresponding GIT condition \cite{Paul-Tian}.
	The Fermat cubic provides a smooth K-stable member \cite[Corollary~4.17]{Zhuang-equivariant}.
	By \cref{thm:stage-I}, every K-semistable $\Q$-Gorenstein limit of smooth cubic $n$-folds is again a cubic hypersurface.
	The moduli continuity argument of Liu--Xu \cite[proof of Theorem~1.1]{Liu-Xu-cubic-threefolds} gives the converse and the asserted isomorphism of stacks.
\end{proof}

\begin{corollary}\label{cor:smooth-cubic-stability}
	Every smooth complex cubic hypersurface of dimension at least two is K-stable.
\end{corollary}

\begin{proof}
	Smooth cubic hypersurfaces are GIT-stable \cite[Proposition~4.2]{Mumford-Fogarty-Kirwan}, so \cref{thm:main} gives the assertion in dimensions $n \geq 5$.
	The cases $n = 2,3,4$ follow from the results of Odaka--Spotti--Sun \cite{OSS}, Liu--Xu \cite{Liu-Xu-cubic-threefolds}, and Liu \cite{Liu-cubic-fourfolds}, respectively.
\end{proof}

\subsection{Proof strategy}
\label{subsec:strategy}

\begin{theorem}\label{thm:stage-I}
	Let $X$ be a K-semistable $\Q$-Gorenstein limit of smooth cubic $n$-folds, where $n \geq 5$.
	Then its limiting hyperplane divisor $L$ is Cartier, and the map $ \phi_{|L|} $ induced by $|L|$ embeds $X$ as a cubic hypersurface in $\PP^{n+1}$.
\end{theorem}

\begin{proof}
	Let $L$ be the ample divisor supplied by \cref{thm:cubic-limits}, so that $-K_X \sim (n-1)L$ and $L^n = 3$.
	The image of $\phi_{|L|}$ has dimension $n$ by \cref{prop:full-dim-image}, so $L$ is Cartier by \cref{prop:full-image-cartier}.
	Fujita's classification of del Pezzo varieties then gives the embedding \cite[Corollary~1.7]{Fujita-del-Pezzo}.
\end{proof}

At a closed non-Cartier point of $L$, the index-one cover associated to $L$ is a singular Gorenstein klt germ with normalized volume greater than $2(n-1)^n$ (\cref{lem:noncartier-cover}).
This exceeds Liu's upper bound for singular local complete intersections \cite[Theorem~1.3]{Liu-cubic-fourfolds}.
The local criterion in \cref{lem:lci-lct-bound} forces the cover to be a local complete intersection if the base ideal of the lifted linear system has log canonical threshold greater than $n-2$.

Let $W$ be the closure of the image of $\phi_{|L|}$.
When $\dim W = n$, the degree of $W$ bounds the mixed intersections on the normalized graph from below.
If the base locus is nonempty, the resulting threshold estimate bounds its codimension from below.
Substituting this codimension bound raises the threshold above $n-2$ (\cref{sec:full-image}).

Suppose now that $m=\dim W<n$.
 The moving part has zero top self-intersection, so the preceding degree estimate no longer applies. 
 Tangent-sheaf semistability applied to the foliation along the fibers gives the negativity needed for adjunction on the Stein base. 
 The section count and the alpha bound then restrict Andreatta’s classification to rational normal scrolls. 
 On the projective bundle model, the two ruling classes sum to the hyperplane class. 
 Applying the same slope bound to both rulings exclude all dimensions except $n=5,6,7$. (See \cref{prop:bounded-lower-image-reduction}).

In the remaining cases, we show that the general fiber of the resolved map is a rational curve. 
A positive-dimensional base locus would force the pullback of $L$ to have degree $\frac{1}{2}$ on this curve.
For the unique exceptional divisor meeting the curve, adjunction would then give a log discrepancy below the mixed-intersection bound in \cref{lem:mixed-intersection-threshold}. 
Thus the base locus is finite, and the same estimate puts the log canonical threshold above $n-2$.
The local criterion gives Cartierness of $L$, and Fujita’s theorem gives very ampleness, contradicting the assumption on the image.(\cref{sec:lower-image}).

\subsection{Related earlier work}

The volume bound underlying the preceding discussion is the following conjecture of Spotti--Sun \cite[Conjecture~5.5]{ODP}; see also \cite[Conjecture~4.5]{Liu-Xu-cubic-threefolds}.

\begin{conjecture}[ODP gap conjecture] \label{conj:odp-gap}
	Let $x\in Z$ be a non-smooth $n$-dimensional klt singularity, with $n\geq2$.
	Then
	\[
		\volhat(x,Z)\leq2(n-1)^n,
	\]
	with equality if and only if $x\in Z$ is an ordinary double point.
\end{conjecture}

Together with the finite degree formula for normalized volumes, this conjecture would give the required Cartierness of the limiting hyperplane divisor directly.
OpenAI has proposed a proof of Conjecture~\ref{conj:odp-gap} in all dimensions \cite{OpenAI-ODP-fourfold,OpenAI-ODP}.
Liu--Xu \cite{Liu-Xu-cubic-threefolds} proved the three-dimensional case using the classification of canonical and terminal singularities.
Liu \cite{Liu-cubic-fourfolds} proved the local complete intersection case in arbitrary dimension and treated cubic fourfolds by studying the restriction of the limiting linear systems to degree-six K3 surfaces.
Our local criterion follows the threshold strategy suggested in \cite[Remark~4.8]{Liu-cubic-fourfolds}, where higher codimensional alpha invariants and surface sections are used to control the index $1$ covers.
We obtain the local complete intersection conclusion from \cite[Lemma~5.7]{Singh-Takagi-Varbaro}.

Our global threshold estimate comes from the fivefold argument of Kong--Shen--Zhao--Zhao \cite{Kong-Shen-Zhao-Zhao-cubic-fivefolds}.
Both arguments use the established volume bound for local complete intersections, rather than the general ODP gap conjecture.

\subsection{An independent proof and a generalization}
\label{subsec:independent-proof}

Kong--Shen--Zhao--Zhao \cite[Theorem~1.1]{Kong-Shen-Zhao-Zhao-cubic-hypersurfaces} give an independent proof of the K-moduli--GIT comparison for cubic hypersurfaces in arbitrary dimension.
Both proofs use Liu's reduction to the Cartierness of the limiting hyperplane divisor.

When $\dim W = n$, they separate cubic and quadric images and analyze the quadrics by rank \cite[proof of Theorem~4.1]{Kong-Shen-Zhao-Zhao-cubic-hypersurfaces}.
We use the codimension of the base locus to obtain a uniform lower bound for the log canonical threshold on the index-one covers.
For $2 \leq \dim W < n$, they use vanishing and Koszul resolutions to show that $W$ has minimal degree \cite[Proposition~3.4]{Kong-Shen-Zhao-Zhao-cubic-hypersurfaces}.
We exclude lower-dimensional images using adjunction on the Stein base and the geometry of a general fiber.

Adding a weighted hyperplane section gives a more general moduli problem.
Si--Zhang--Zhou \cite{Si-Zhang-Zhou-cubic-hyperplane-pairs} formulated the K-moduli--GIT comparison for these pairs.
Let $X$ be a cubic $n$-fold with $n \geq 2$, and let $D$ be a hyperplane section.
For rational $0 < t < 2/n$, they show that every K-semistable limit of such pairs is again a cubic–hyperplane pair. 
They also identify the good moduli space of the reduced smoothable K-moduli component of pairs $(X,tD)$ with the GIT quotient at slope $s(t) = \frac{3(n+1)t}{2n+4+t}$.
The semistable, polystable, and stable loci agree, and the identifications commute with wall crossing.
For $n \geq 3$, the comparison also holds at $t = 2/n$.

\subsection*{Acknowledgements}

We first thank Junyan Zhao for the discussions about the main ideas of their proof and for kindly sharing their draft \cite{Kong-Shen-Zhao-Zhao-cubic-hypersurfaces}.
We also thank Zhiyuan Li and Chenyang Xu for their interest in our results and their valuable mathematical advice.
L.P thank Chuyu Zhou for inviting him to visit Xiamen University and for sharing the draft on cubic--hyperplane pairs \cite{Si-Zhang-Zhou-cubic-hyperplane-pairs}.
L.P was partially supported by the CAS Project for Young Scientists in Basic Research (Grant No.YSBR-033), and Postdoctoral Fellowship Program of CPSF (Grant No.GZC20261659).
H.W was partially supported by China Postdoctoral Science Foundation (Grant No.043210052).

\subsection*{AI usage disclosure}

\label{subsec:AI-disclosure}
We are sharing this manuscript with everyone interested in the problem, in the spirit of sharing mathematical knowledge.
\[
	\text{\emph{
			We do not claim any traditional academic credit for this paper.}}\]
\[
	\text{\emph{This paper is not intended for publication in a traditional mathematics journal.}}\]

LLM(AI) contributed substantive mathematical suggestions, including proof strategies and arguments used in preparing the manuscript.
Its principal contribution concerned the lower-dimensional image case as they suggested combining the foliations defined by the rational maps with tangent-sheaf semistability and adjunction theory.

The authors reviewed the resulting arguments, removed unnecessary detail, and revised the exposition.
In particular, the analysis in \cref{sec:full-image} and \cref{sec:lower-image} was simplified.
We take full responsibility for the mathematical content and its presentation, including any remaining errors or omissions.

This work owes a fundamental intellectual debt to Yuchen Liu and Chenyang Xu.
Liu's study of cubic degenerations through the limiting hyperplane divisor and its linear system is the starting point of our argument.
Xu's work with his collaborators on algebraic K-stability and K-moduli supplies the framework in which it is carried out.
The use of AI tools does not alter the attribution of these ideas or the results of other researchers used below.

\section{Preliminaries}
\label{sec:preliminaries}

We work over $\C$.
Unless otherwise specified, all varieties are assumed to be normal and projective.
However, we do not assume normality of the image of a rational map.
We always write $X$ for a $\Q$-Gorenstein limit of smooth cubic hypersurfaces.

\subsection{K-stability and K-moduli}
\label{subsec:valuative}

\begin{definition}[$\Q$-Fano varieties]\label{def:qfano}
	A \emph{$\Q$-Fano variety} is a normal projective klt variety $Z$ with ample $\Q$-Cartier anticanonical divisor $-K_Z$.
	A \emph{log Fano pair} $(Z,\Delta)$ consists of a projective klt pair with effective $\Q$-divisor $\Delta$ and ample $-(K_Z+\Delta)$.
\end{definition}

\begin{definition}[Log canonical thresholds]
	\label{def:lct-systems}
	Let $Z$ be a $\Q$-Gorenstein klt variety.
	For a nonzero coherent ideal $\mathfrak{a}\subset\OK_Z$, set
	\[
		\lct(Z;\mathfrak{a})
		:=\inf_{E}
		\frac{A_Z(E)}{\ord_E(\mathfrak{a})}.
	\]
	We set $\frac{A_Z(E)}{\ord_E(\mathfrak{a})}=+\infty$ if $\ord_E(\mathfrak{a})=0$.
	The log canonical threshold $\lct_z(Z;\mathfrak{a})$ is defined by divisors whose centers contain $z$.
	For an ideal $I \subset \OK_{Z,z}$, we also use the notation $\lct_z(Z;I)$. This denotes the local threshold of any coherent ideal sheaf representing $I$ near $z$.
	For an effective $\Q$-Cartier divisor $D$, or a linear system $\mathcal{V}$ of such divisors, we replace the denominator by $\ord_E(D)$, or $\ord_E(\mathcal{V})=\min_{D'\in\mathcal{V}}\ord_E(D')$ in the denominator.
	When $\mathcal{V}$ is a subsystem of a Cartier linear system, its threshold is that of $\mathfrak{b}(\mathcal{V})$.
	It satisfied the following properties:
	\[
		\lct(Z;\lambda\mathcal{V})
		=\lambda^{-1}\lct(Z;\mathcal{V})\qquad(\lambda>0).
	\]

\end{definition}

\begin{definition}[Valuative invariants]\label{def:invariants}
	Let $Z$ be an $n$-dimensional $\Q$-Fano variety.
	Let $L$ be an ample $\Q$-Cartier divisor, and let $E$ be a prime divisor on a smooth projective birational model $\pi:Y\to Z$.
	Define the $S$-invariant:
	\[
		S_L(E)=\frac{1}{L^n}\int_0^{\infty}
		\vol\bigl(\pi^*L-tE\bigr)\,dt.
	\]
	Changing variables in the integral gives $S_{cL}(E)=cS_L(E)$ for
	$c\in\Q_{>0}$.
	We write $S_Z(E)=S_{-K_Z}(E)$ and define the K-stability threshold by
	\[
		\delta(Z)=\inf_E\frac{A_Z(E)}{S_Z(E)},
	\]
	where $E$ runs over all prime divisors over $Z$.
	For a log Fano pair $(Z,\Delta)$, one defines $\delta(Z,\Delta)$ by using $A_{Z,\Delta}(E)$ and $S_{-(K_Z+\Delta)}(E)$.
\end{definition}

\begin{theorem}[{\cite{Fujita-valuative,Li-Kss-equivariant}}]
	\label{ext:valuative}
	A $\Q$-Fano variety $Z$ is K-semistable if and only if $\delta(Z)\geq1$, equivalently if and only if
	\[
		A_Z(E)\geq S_{-K_Z}(E)
	\]
	for every prime divisor $E$ over $Z$.
\end{theorem}

\begin{definition}[Alpha invariants]\label{def:alpha-invariants}
	For an ample $\Q$-Cartier divisor $L$ on an $n$-dimensional projective klt variety $Z$, we define
	\[
		\alpha(Z;L)
		=\inf\{\lct(Z;D)\mid D\geq0,\ D\sim_{\Q}L\},
	\]
	where $D$ ranges over effective $\Q$-divisors.
	For $s \in \Q_{>0}$, one has $\alpha(Z;sL) = s^{-1}\alpha(Z;L)$.
	If $Z$ is $\Q$-Fano, we write $\alpha(Z)=\alpha(Z;-K_Z)$.
\end{definition}

\begin{theorem}[{\cite[Theorem~1.1, $k = 1$]{Zhu-higher-alpha}}]
	\label{ext:higher-alpha}
	If $Z$ is an $n$-dimensional K-semistable $\Q$-Fano variety, then $\alpha(Z) \geq 1/(n+1)$.
	In particular, if $-K_Z\sim_{\Q}(n-1)L$, then
	$\alpha(Z;L)\geq(n-1)/(n+1)$.
\end{theorem}

The comparison with GIT uses three facts from K-moduli theory: existence of polystable replacements, properness of good moduli spaces, and openness of K-semistability.
Furthermore, people prove the following K-moduli theorem:

\begin{theorem}[{\cite[Theorem~1.3 and Remark~1.5]{Liu-Xu-Zhuang-Finitegeneration}}]
	Fix $n\in\mathbb{N}_{>0}$ and $v\in\mathbb{Q}_{>0}$.
	The moduli functor which sends any scheme $S$ to
	\[
		\left\{
		\mathcal{X}/S\Bigg|
		\begin{array}{cc}
		\mathcal{X}\rightarrow S\text{ is a proper flat family satisfying Koll\'{a}r condition,}     & \\
		\text{each fiber $\mathcal{X}_s$ is an $n$-dimensional K-semistable $\mathbb{Q}$-Fano variety} & \\
		\text{ of volume $v$ } 
		\end{array}\right\} \Bigg/ \sim_{\rm iso} 
	\]
	can be represented by an Artin stack of finite type over $\mathbb{C}$, which is  denoted by $\mathcal{M}^{K}_{n,v}$.
	Moreover, it admits a projective good moduli space $M^{K}_{n,v}$ parametrizing all K-polystable $\mathbb{Q}$-Fano varieties with fixed dimension $n$ and volume $v$.
\end{theorem}

\begin{remark}[Smoothable families]\label{rem:kmoduli-existence-source}
	In the $\Q$-Gorenstein smoothable setting, proper good moduli spaces were constructed in \cite[Theorem~1.3]{Li-Wang-Xu-proper-moduli}.
	For smoothable log pairs with boundary proportional to the anticanonical divisor, existence of K-semistable fillings is established in \cite[Theorem~3.19]{Ascher-DeVleming-Liu}.
	Uniqueness of K-polystable limits in families is proved in \cite{Blum-Xu-uniqueness}.
\end{remark}

\subsection{Singularities and normalized volumes}
\label{subsec:volume-bounds}

\begin{definition}[Normalized volumes]\label{def:normalized-volume}
	Let $x\in Z$ be an $n$-dimensional klt germ and let $E$ be a prime divisor over $Z$ centered at $x$.
	Its valuation ideals and volume are defined by:
	\[
		\mathfrak{a}_j(E)
		=\{f\in\OK_{Z,x}\mid\ord_E(f)\geq j\},\qquad
		\vol_{Z,x}(E)
		=\limsup_{j\to\infty}
		\frac{\ell\bigl(\OK_{Z,x}/\mathfrak{a}_j(E)\bigr)}{j^n/n!},
	\]
	where $\ell$ denotes length.
	Define
	\[
		\volhat_{Z,x}(E)=A_Z(E)^n\vol_{Z,x}(E),\qquad
		\volhat(x,Z)=\inf_{c_Z(E)=x}\volhat_{Z,x}(E).
	\]
	The latter is called the local volume of  the germ.
\end{definition}

\begin{theorem}
	\label{ext:volume-bounds}
	The following normalized-volume estimates hold.
	\begin{enumerate}
		\item \textup{\cite[Theorem~2.14]{Liu-cubic-fourfolds}.}
		      If $Z$ is an $n$-dimensional K-semistable $\Q$-Fano variety, then every closed point $x\in Z$ satisfies
		      \[
			      (-K_Z)^n\leq\Bigl(1+\frac1n\Bigr)^n\volhat(x,Z).
		      \]
		\item \textup{\cite[Theorem~2.13]{Liu-cubic-fourfolds}.} If $\tau:(y\in Y)\to(x\in Z)$ is a finite quasi-\'etale Galois morphism between klt singularities, then
		      $\volhat(y,Y)=\deg(\tau)\volhat(x,Z)$.
		\item \textup{\cite[Theorem~2.12]{Liu-cubic-fourfolds}.} If $y\in Y$ is a non-smooth local complete intersection klt singularity of dimension $n$, then
		      \[
			      \volhat(y,Y)\leq2(n-1)^n,
		      \]
		      with equality if and only if it is an ordinary double point.
	\end{enumerate}
\end{theorem}

\subsection{Polarized models and foliations}
\label{subsec:polarized-foliations}

\begin{definition}[Polarized birational models]\label{def:polarized-models}
	A quasi-polarized pair $(B,A)$ consists of a normal projective variety and a nef and big Cartier divisor.
	Two such pairs are said to be polarized birationally equivalent if their polarizations have the same pullback to a common resolution; see \cite[Definition~3.1]{Andreatta}.
	This identifies their section spaces and images of complete linear-series maps.

	For a terminal $\Q$-factorial pair of dimension $m\geq2$, a \textit{zero reduction} is the last pair in a $(K_B+\Delta)$-minimal model program with scaling, where $(B,\Delta)$ is klt and
	$\Delta\sim_{\Q}(m-1)A$.
	The polarization descends through each birational step, so a zero reduction is polarized birationally equivalent to $(B,A)$; see \cite[Proposition~4.2, Definition~4.4, and Remark~4.5(i)]{Andreatta}.
\end{definition}

\begin{definition}[Slopes]\label{def:slopes}
	Let $L$ be an ample $\Q$-Cartier divisor on an $n$-dimensional normal projective variety $Z$.
	The slope with respect to $L$ of a torsion-free sheaf $\mathcal{E}$ is defined by
	\[
		\mu_L(\mathcal{E})
		=\frac{c_1(\mathcal{E})\cdot L^{n-1}}{\rk(\mathcal{E})}.
	\]
	Here $c_1(\mathcal{E})$ is the Weil divisor class of the reflexive determinant
	$\det\mathcal{E}:=(\bigwedge^{\rk(\mathcal{E})}\mathcal{E})^{**}$;
	its intersection with $L^{n-1}$ is defined even if that class is not $\Q$-Cartier.
	The sheaf is said to be slope semistable if $\mu_L(\mathcal{F})\leq\mu_L(\mathcal{E})$ for every coherent subsheaf with $0<\rk(\mathcal{F})<\rk(\mathcal{E})$.
	Scaling $L$ by a positive rational number does not change semistability.
\end{definition}

\begin{definition}[Foliations]\label{def:foliation}
	Let $Z$ be a normal complex projective $n$-fold with K\"ahler differentials $\Omega_Z^1$ and reflexive tangent sheaf $T_Z = (\Omega_Z^1)^*$.
	A subsheaf $\mathcal{E} \subset T_Z$ is \textit{saturated} if $T_Z/\mathcal{E}$ is torsion-free. Its saturation is the smallest saturated subsheaf containing $\mathcal{E}$.
	A \textit{foliation} on $Z$ is a saturated subsheaf $\mathcal{F} \subset T_Z$ closed under the Lie bracket.
	Its leaves are the maximal connected integral submanifolds on the locus where $Z$ and $\mathcal{F}$ are regular.

	Let $\phi\colon Z \dashrightarrow B$ be a dominant rational map to a smooth projective $b$-fold, with $0 < b < n$.
	Its \textit{kernel foliation} is the saturated subsheaf $\mathcal{F} \subset T_Z$ agreeing with $\ker(d\phi)$ wherever $\phi$ is defined and smooth.
	Take any resolution $p\colon Y \to Z$  and $q \colon  Y \to B$ of $\phi$. 
	Then $\mathcal{F}$ is the saturation of image of $p_*T_{Y/B}$ in $T_Z$.
\end{definition}

\begin{lemma}\label{lem:kernel-foliation-bound}
	With the notation of \cref{def:foliation}, $c_1(\mathcal{F}) = -K_Z + p_*q^*K_B + R$ in $\operatorname{Cl}(Z)$ for some effective divisor $R$ on $Z$.
	If $T_Z$ is slope semistable with respect to an ample $\Q$-Cartier divisor $L$, then
	\[
		q^*(-K_B) \cdot (p^*L)^{n-1}
		\geq \frac{b}{n}(-K_Z) \cdot L^{n-1}.
	\]
\end{lemma}

\begin{proof}
	Over the codimension-one points of $Z$, the map $p$ is an isomorphism and $dq$ embeds $T_Z/\mathcal{F}$ into $q^*T_B$ with full rank.
	The determinant of this inclusion vanishes along an effective divisor $R$, giving the first formula as in \cite[Section~6, equation~(8)]{Horing-Peternell}.
	The inequality follows from $\rk(\mathcal{F}) = n-b$, $\mu_L(\mathcal{F}) \leq \mu_L(T_Z)$ and the projection formula.
\end{proof}

\section{A local criterion for Cartierness}
\label{sec:thresholds_at_limits}

In this section, we first recall the structure theorem for cubic limits and the associated limiting hyperplane divisor $L$.
At a closed non-Cartier point $x$ of $L$, the volume estimate in \cref{lem:noncartier-cover} shows that the index-one cover is not a local complete intersection.
Combining this with \cref{lem:lci-lct-bound} gives a local Cartier criterion in terms of $\lct_x(X;|L|)$ (\cref{cor:local-cartier-criterion}).

\subsection{Cubic limits and index-one covers}
\label{subsec:cubic-limits}

We start with Liu's structure theorem.
It supplies the limiting hyperplane class and reduces the possible failure of Cartierness to the case of local index two.

\begin{theorem}[{\cite[Theorem~3.1]{Liu-cubic-fourfolds}}]
	\label{thm:cubic-limits}
	Let $n\geq4$ and let $\mathcal{X}\to B$ be a K-semistable $\Q$-Fano family over a smooth pointed curve $0\in B$, whose restriction over $B\setminus\{0\}$ is a smooth family of cubic $n$-folds.
	Put
	$X=\mathcal{X}_0$.
	There is a $\Q$-Cartier integral Weil divisor class $L$ on $X$ with the following properties.
	\begin{enumerate}
		\item $\OK_X(mL)$ is Cohen--Macaulay for every integer $m$.
		\item $-K_X\sim(n-1)L$ and $L^n=3$.
		\item $h^i(X,\OK_X(mL))=h^i(\mathcal{X}_b,\OK_{\mathcal{X}_b}(m))$ for every integer $m$, every $i\geq0$ and every $b\neq0$.
		      In particular $h^0(X,L)=n+2$ and the intermediate cohomology vanishes for all $m$.
		\item Every $\Q$-Cartier Weil divisor $D$ on $X$ has $2D$ Cartier.
	\end{enumerate}
\end{theorem}

\begin{definition}[The cubic limit pair]\label{def:cubic-limit-pair}
	Throughout the paper, $(X,L)$ denotes a K-semistable $\Q$-Gorenstein limit of smooth cubic $n$-folds with $n \geq 5$, together with the ample integral $\Q$-Cartier divisor $L$ supplied by \cref{thm:cubic-limits}, unless other hypotheses are stated.
	In particular, the non-Cartier locus of $L$ is contained in the geometric base locus $\mathcal{B}(|L|)$.
\end{definition}

The next lemma gives a sufficient condition for Cartierness of $L$ at a closed point in terms of the singularities of the index-one cover.
The volume comparison is immediate once we know that the cover is singular.

\begin{lemma}
	\label{lem:noncartier-cover}
	In the setting of \cref{thm:cubic-limits}, suppose $L$ is not Cartier at a closed point $x\in X$, and let $\tau:(y\in Y)\to(x\in X)$ be its index-one cover.
	Then $\deg(\tau) = 2$, $(Y,y)$ is a singular Gorenstein klt germ, and
	\[
		\volhat(y,Y)\geq6\Bigl(\frac{n(n-1)}{n+1}\Bigr)^n>2(n-1)^n.
	\]
	In particular, $Y$ is not a local complete intersection at $y$.
\end{lemma}

\begin{proof}
	By \cref{thm:cubic-limits}(4), the local Cartier index of $L$ at $x$ is two.
	Thus $\tau$ is a quasi-\'etale double cover with a unique point $y$ over $x$.
	Since $K_X \sim -(n-1)L$ and $\tau^{[*]}L$ is Cartier, $K_Y = \tau^*K_X$ is Cartier at $y$.
	The quasi-\'etale cover is klt and hence Cohen--Macaulay, so it is Gorenstein at $y$.

	We next show that $Y$ is not smooth at $y$.
	Suppose, for a contradiction, that $Y$ is smooth at $y$.
	After a finite base change if necessary, the construction in the proof of \cite[Theorem~3.1]{Liu-cubic-fourfolds} provides a Weil divisor $\mathcal{L}$ on $\mathcal{X}$ extending $L$, such that $\OK_{\mathcal{X}}(j\mathcal{L})$ is Cohen--Macaulay and restricts to $\OK_X(jL)$ for every $j\in\Z$.
	In particular,
	\[
		\OK_{\mathcal{X}}(2\mathcal{L})/
		t\OK_{\mathcal{X}}(2\mathcal{L})
		\cong \OK_X(2L),
	\]
	where $t$ is a local parameter of $B$ at $0$.
	The right-hand side is invertible, so Nakayama's lemma shows that the rank-one torsion-free module $\OK_{\mathcal{X}}(2\mathcal{L})_x$ is cyclic, hence free over the local domain $\OK_{\mathcal{X},x}$.
	Thus $2\mathcal{L}$ is Cartier near $x$.

	After shrinking around $x$, we can choose a trivialization of $\OK_{\mathcal{X}}(2\mathcal{L})$ compatible with the one defining $\tau$.
	The associated cyclic algebra
	\[
		\mathcal{A}
		=
		\OK_{\mathcal{X}}
		\oplus
		\OK_{\mathcal{X}}(\mathcal{L})
	\]
	defines a finite double cover
	$\mathcal{Y}\to\mathcal{X}$.
	Since both summands of $\mathcal{A}$ are flat over $B$, the morphism $\mathcal{Y}\to B$ is flat.
	The restriction isomorphisms identify its central fiber with the original index-one cover $\tau:Y\to X$.
	Over $t\neq 0$, the cover is \'etale because $\mathcal{L}$ is Cartier there.

	Since $\mathcal{Y} \to B$ is flat and its central fiber $Y$ is smooth at $y$, this morphism is smooth near $y$.
	Let $\iota$ be the deck involution, and shrink to an $\iota$-invariant neighborhood of $y$.
	The involution acts over $B$, and its order is invertible, so its fixed locus $F$ is smooth over $B$ by \cite[Proposition~3.4]{Edixhoven-neron-models}.
	Since $y \in F$, the image of $F \to B$ contains an open neighborhood of $0$.
	Thus $F$ meets a noncentral fiber, contradicting the freeness of the deck action over $B \setminus \{0\}$.

	Now we can apply the local-to-global inequality and the finite degree formula in \cref{ext:volume-bounds} (1)--(2).
	Since $(-K_X)^n=3(n-1)^n$, this gives the asserted lower bound.
	Its ratio to $2(n-1)^n$ is
	$3\bigl(n/(n+1)\bigr)^n>1$.
	Estimate~(3) in \cref{ext:volume-bounds} therefore excludes the l.c.i. property at $y$.

\end{proof}

\begin{remark}
	If we assume the ODP gap conjecture \ref{conj:odp-gap}, then we directly conclude that $L$ is Cartier at every closed point of $X$.
	Both \cref{thm:stage-I} and \cref{thm:main} then follow immediately.
\end{remark}

\subsection{Local thresholds on index-one covers}
\label{subsec:local-threshold-bounds}

The local rings below are essentially of finite type over a characteristic-zero field.
At a nonclosed point, lifts of a separating transcendence basis of the residue field realize the same local ring at a closed point of a generic fiber. A log resolution restricts to this fiber with the same discrepancies and vanishing orders along divisors meeting it.
Extension to an algebraic closure preserves klt singularities and thresholds \cite[Lemma~2.14]{deFernex-Ein-Mustata-bounded}, as well as Gorensteinness, local dimension, and embedding dimension.

\begin{lemma}\label{lem:lci-lct-bound}
	Let $(Z,z)$ be a Gorenstein klt local germ of dimension $c\geq2$ over a characteristic-zero field, and let $I\subsetneq\OK_{Z,z}$ be a nonzero ideal.
	Assume
	\[
		\lct_z(Z;I)>c-2.
	\]
	Then $Z$ is a local complete intersection at $z$.
	More precisely, its completed local ring is regular or a hypersurface.
\end{lemma}

\begin{proof}
	Write $R=\OK_{Z,z}$, with maximal ideal $\mathfrak{m}$ and residue field $\kappa$.
	Since $I\subseteq\mathfrak{m}$,
	\[
		t:=\lct_z(Z;\mathfrak{m})\geq\lct_z(Z;I)>c-2.
	\]
	Choose a minimal reduction $Q=(q_1,\dots,q_c)$ of $\mathfrak{m}$.
	This exists because $\kappa$ is infinite \cite{Northcott-Rees}.
	By \cite[Lemma~5.7(2), (3)]{Singh-Takagi-Varbaro},
	\[
		t\leq c,\qquad
		\mathfrak{m}^{\,c+1-\lceil t\rceil}\subseteq Q.
	\]
	The exponent is at most two, so $\mathfrak{m}^2\subseteq Q$.

	Since $R$ is Gorenstein and $Q$ is a parameter ideal, $A=R/Q$ is Artinian Gorenstein.
	Its maximal ideal $\mathfrak{n}=\mathfrak{m}/Q$ satisfies $\mathfrak{n}^2=0$, hence
	\[
		\mathfrak{n}\subseteq\operatorname{Soc}(A)\coloneqq \operatorname{ann}_A(\mathfrak{n}),\qquad
		\dim_{\kappa}\operatorname{Soc}(A)=1.
	\]
	Lifting at most one generator of $\mathfrak{n}$ and adjoining $q_1,\dots,q_c$ therefore gives
	\[
		\operatorname{edim}(R)\leq c+1.
	\]
	Since $R$ is Cohen--Macaulay, a minimal Cohen presentation of $\widehat{R}$ has kernel either zero or generated by one nonzerodivisor.
	Thus $\widehat{R}$ is regular or a hypersurface, and $Z$ is a local complete intersection at $z$.
\end{proof}

Let $L$ be an integral $\Q$-Cartier Weil divisor on a normal projective klt variety $Z$, and put $V=H^0(Z,L)\neq0$.
If $L$ has local Cartier index $r$ at $z\in Z$, denote by
\[
	\tau:(\widetilde z\in\widetilde Z)\longrightarrow(z\in Z)
\]
the index $1$ cover associated to $L$, of degree $r$.
All such covers are understood as covers of germs.
Put $\widetilde L=\tau^{[*]}L$, which is Cartier.

\begin{definition}\label{def:thresholds-on-covers}
	The pullbacks of the sections in the fixed space $V$ define the coherent ideal sheaf
	\[
		\mathfrak{I}
		=\im\bigl(V\otimes\OK_{\widetilde Z}(-\widetilde L)
		\longrightarrow\OK_{\widetilde Z}\bigr),\qquad
		I_{\widetilde z}\coloneqq \mathfrak{I}_{\widetilde z}
		\subset\OK_{\widetilde Z,\widetilde z}.
	\]
	We call $\mathfrak{I}$ the \emph{base ideal of the lifted linear system}, or simply the \emph{base ideal}.
	In local statements, the same term refers to its stalk $I_{\widetilde z}$.
\end{definition}

After trivializing $\widetilde L$ near $\widetilde z$, the base ideal is generated by the functions representing the pullbacks of all global sections of $L$.
Where $L$ is Cartier, we take $r=1$ and the identity cover.

Let $w$ be a normalized divisorial valuation over $\widetilde Z$ whose center contains $\widetilde z$.
Write $w|_{\C(Z)} = e\ord_G$, where $G$ is a prime divisor over $Z$ and $e = e(w/G)$ is the ramification index.
Since the cover is cyclic of degree $r$, we have $e\mid r$.
Quasi-\'etaleness and crepancy give
\[
	\begin{aligned}
		A_{\widetilde Z}(w) & = e A_Z(G),     \\
		w(I_{\widetilde z}) & = e\ord_G(|L|),
	\end{aligned}
\]
where $\ord_G(|L|) = \min_{D\in|L|}\ord_G(D)$.
Conversely, by finiteness and Galois conjugacy, every $\ord_G$ whose center contains $z$ admits a divisorial extension to $\C(\widetilde{Z})$ with center containing $\widetilde{z}$.
Thus the ratio of log discrepancy to vanishing order is preserved, and
\[
	\lct_{\widetilde z}(\widetilde Z;I_{\widetilde z})
	= \lct_z(Z;|L|).
\]

\begin{corollary}[Local Cartier criterion]\label{cor:local-cartier-criterion}
	Let $(X,L)$ be a cubic limit as in \cref{thm:cubic-limits}, with $\dim X = n \geq 4$, and let $x \in X$ be a closed point.
	If $\lct_x(X;|L|) > n-2$, then $L$ is Cartier at $x$.
\end{corollary}

\begin{proof}
	Suppose that $L$ is not Cartier at $x$.
	By \cref{lem:noncartier-cover}, its index-one cover $(\widetilde{X},\widetilde{x})$ is Gorenstein klt and is not a local complete intersection.
	Every member of $|L|$ passes through $x$, so the lifted sections vanish at $\widetilde{x}$ and $I_{\widetilde{x}}$ is a proper ideal.
	The threshold comparison above gives $\lct_{\widetilde{x}}(\widetilde{X};I_{\widetilde{x}}) > n-2$, so \cref{lem:lci-lct-bound} forces the cover to be a local complete intersection at $\widetilde{x}$, a contradiction.
\end{proof}

\section{Linear systems and mixed intersections}
\label{sec:full-image}
In this section, we prove the Cartierness of the limiting hyperplane divisor $L$ in \cref{def:cubic-limit-pair} when the image of the linear system has full dimension.

\begin{proposition}\label{prop:full-image-cartier}
	Let $(X,L)$ be a K-semistable cubic limit as in \cref{thm:cubic-limits}, with $\dim X = n \geq 4$.
	Assume the image of $\phi_{|L|}$ has dimension $n$.
	Then $L$ is Cartier.
\end{proposition}

Throughout this section, let $Z$ be a normal projective variety of dimension $n\geq2$, and let $L$ be an ample integral $\Q$-Cartier Weil divisor with $H^0(Z,L)\neq0$.

\subsection{The normalized graph}
\label{subsec:product-ideal-model}
Set
\[
	V=H^0(Z,L),\qquad N=\dim V-1,\qquad
	\phi=\phi_{|L|},\qquad W=\overline{\phi(Z)}\subset\PP^N,
\]
and write $m=\dim W$.
Write $\mathcal{B} = \bigcap_{D\in|L|}\Supp D$ for the geometric base locus of $|L|$.
Choose a positive integer $\ell$ such that $\ell L$ is Cartier.
Since we use the complete linear system, $W$ is integral and nondegenerate; we make no assumption on its normality.

Let
\[
	V_\ell=\im\bigl(\Sym^\ell V\longrightarrow H^0(Z,\ell L)\bigr),\qquad
	\mathfrak{b}_\ell
	=\im\bigl(V_\ell\otimes\OK_Z(-\ell L)\longrightarrow\OK_Z\bigr).
\]
The zero set of $\mathfrak{b}_\ell$ is $\mathcal{B}$, since $V_\ell$ is spanned by products of $\ell$ sections of $V$ and contains every $\ell$-th power.
On an index-one cover $\tau$, the identification $\tau^*\OK_Z(\ell L) \simeq \OK_{\widetilde{Z}}(\ell\widetilde{L})$ takes these generators to products of $\ell$ lifted sections.
Thus $\mathfrak{b}_\ell\OK_{\widetilde{Z}} = \mathfrak{I}^{\ell}$ for the ideal of \cref{def:thresholds-on-covers}.

Take $Y$ to be the normalization of the blowup of $\mathfrak{b}_\ell$. 
Then we can identify $Y$ with the normalized graph of $\phi$, independent of $\ell$, with morphisms $p:Y\to Z$ and $f:Y\to W$ fitting into the diagram
\[
	\begin{tikzcd}
		& Y \arrow[dl, "p"'] \arrow[dr, "f"] & \\
		Z \arrow[rr, dashed, "\phi"'] & & W\subset\PP^N.
	\end{tikzcd}
\]
Let $A_W=\OK_W(1)$.
Put $P=p^*L$ and $M=f^*A_W$.
Both $P$ and $M$ are nef, and $M$ is Cartier, globally generated, and $p$-ample.
Moreover,
\[
	P=M+F,\qquad F\geq0,\qquad \mathfrak{b}_\ell\OK_{Y}=\OK_{Y}(-\ell F) \text{~Cartier}.
\]
The support satisfies $\Supp F=p^{-1}(\mathcal{B})$ set-theoretically, and $F$ is the fixed $\Q$-divisor of $P$, i.e.,
\[
	H^0(Z,L)\cong H^0(Y,M).
\]

\subsection{Mixed intersections and threshold estimates}
\label{subsec:mixed-thresholds}
Define the mixed intersections on the normalized graph $Y$ by
\[
	a_i=P^{n-i}M^i\qquad(0\leq i\leq n).
\]
They are unchanged by pullback to any higher birational model.
On such a model, we continue to write $P,M,F$ for their pullbacks; when comparing two models, we write the pullbacks explicitly.
For a prime divisor $G$ over $Z$, let $b_G$ be the minimal vanishing order along $G$ of a member of $|L|$, i.e.
\[
	b_G\coloneqq\ord_G(|L|)=\min_{D\in|L|}\ord_G(D).
\]
On a model carrying $G$, this is the coefficient of $G$ in the pull-back of the fixed divisor $F$.

We first show how the mixed intersections reflect the dimension of the image and the codimension of the base locus.

\begin{lemma}\label{lem:base-codimension-mixed-intersections}
	The mixed intersections satisfy
	\[
		L^n=a_0\geq a_1\geq\cdots\geq a_n\geq0,\qquad
		a_i=0\quad\text{if }i>m.
	\]
	If $\mathcal{B}\neq\varnothing$ and $k=\codim_Z\mathcal{B}$, then
	\[
		a_0=\cdots=a_{k-1}=L^n>a_k.
	\]
\end{lemma}

\begin{proof}
	Since $F=P-M$ is effective and $P,M$ are nef,
	\[
		a_i-a_{i+1}=F\cdot P^{n-i-1}M^i\geq0
		\qquad(0\leq i\leq n-1).
	\]
	Also $a_0=L^n$, and $a_i=0$ for $i>m$ because $M$ is pulled back from $W$.

	Suppose $\mathcal{B}$ has codimension $k$.
	For $0\leq i\leq k-2$, the effective cycle $F\cdot M^i$ has dimension $n-1-i>n-k$, whereas its image under $p$ lies in $\mathcal{B}$.
	Hence $p_*(F\cdot M^i)=0$.
	The projection formula gives
	\[
		a_i-a_{i+1}
		=L^{n-i-1}\cdot p_*(F\cdot M^i)=0,
	\]
	proving the equalities.
	For $k=1$ there is only $a_0=L^n$.

	To prove the strict inequality,   choose an irreducible component $T$ of $\mathcal{B}$ of codimension $k$, and let $\eta_T$ be its generic point.
	Since $\Supp F=p^{-1}(\mathcal{B})$, some prime component $D$ of $F$ dominates $T$.
	Its generic fiber $D_{\eta_T}$ has dimension $k-1$.
	The projection formula gives
	\[
		D\cdot P^{n-k}M^{k-1}
		=(L^{n-k}\cdot T)\,
		\deg\bigl((M|_{D_{\eta_T}})^{k-1}\bigr)>0
	\]
	by the ampleness of $L$ and the $p$-ampleness of $M$.
	All components of $F$ have nonnegative intersection with the nef classes $P^{n-k}M^{k-1}$, and the coefficient of $D$ in $F$ is positive.
	Hence
	\[
		a_{k-1}-a_k=F\cdot P^{n-k}M^{k-1}>0.
	\]
	This is the positivity argument of \cite[proof of Theorem~3.2, Step~2]{Horing-Peternell}.
\end{proof}

Thus the codimension of a nonempty base locus is the index of the first strict drop in the mixed-intersection sequence.
In particular, equality of the first $n$ mixed intersections characterizes finiteness of the base locus.



The segment joining $P$ and $M$ gives a lower bound for the volume integral defining $S_L$ in terms of the mixed intersections.

\begin{lemma}\label{lem:mixed-intersection-threshold}
	For every prime divisor $G$ over $Z$,
	\[
		S_L(G)\geq\frac{b_G}{(n+1)L^n}\sum_{i=0}^n a_i.
	\]
	Assume in addition that $Z$ is a K-semistable $\Q$-Fano variety and $-K_Z\sim_{\Q}\lambda L$ for some
	$\lambda\in\Q_{>0}$.
	Set
	\[
		\theta=\frac{\lambda}{(n+1)L^n}\sum_{i=0}^n a_i.
	\]
	Then
	\[
		A_Z(G)\geq S_{-K_Z}(G)\geq\theta\cdot b_G.
	\]
\end{lemma}

\begin{proof}
	Choose a smooth model $\sigma:Y'\to Y$ carrying $G$, and use $P,M,F$ for the pulled-back divisors.
	Then $F-b_GG\geq0$.
	The volume bound is immediate if $b_G=0$.
	Otherwise, for $0\leq t\leq1$,
	\[
		P-t\cdot b_GG=(1-t)P+tM+t(F-b_GG)\geq(1-t)P+tM.
	\]
	The divisor on the right is nef.
	Monotonicity of volume and integration over $0\leq tb_G\leq b_G$ give
	\[
		S_L(G) \geq\frac{b_G}{L^n}\int_0^1\bigl((1-t)P+tM\bigr)^n\,dt =\frac{b_G}{(n+1)L^n}\sum_{i=0}^n a_i.
	\]
	Here
	$\int_0^1\binom{n}{i}(1-t)^{n-i}t^i\,dt=1/(n+1)$.
	Under the additional K-semistable Fano hypotheses, $S_{-K_Z}(G)=\lambda S_L(G)$, and the valuative criterion proves the discrepancy bound.
\end{proof}

At the generic point of a base component, the codimension of the component bounds the local threshold of $|L|$ from above.

\begin{lemma}\label{lem:codim_lct-intersection-bound}
	Assume that $Z$ is a K-semistable $\Q$-Fano variety and $-K_Z\sim_{\Q}\lambda L$ for some $\lambda\in\Q_{>0}$, and let $\theta$ be as in \cref{lem:mixed-intersection-threshold}.
	If $\mathcal{B}$ is nonempty, then
	\[
		\codim_Z \mathcal{B} \geq \lct(Z;|L|) \geq \theta.
	\]
	For a cubic limit $(X,L)$ as in \cref{thm:cubic-limits}, $L$ is Cartier whenever $\theta > n-2$.
\end{lemma}

\begin{proof}
	The discrepancy bound in \cref{lem:mixed-intersection-threshold} gives $\lct_z(Z;|L|) \geq \theta$ at every scheme point $z \in Z$. Taking the infimum gives the global bound by \cref{def:lct-systems}.

	If $\mathcal{B}$ is nonempty, let $\eta$ be the generic point of a component of codimension $k = \codim_Z \mathcal{B}$.
	Let $(\widetilde{Z},\widetilde{\eta}) \to (Z,\eta)$ be the index-one cover of $L$.
	The cover is finite, so the klt local ring at $\widetilde{\eta}$ has dimension $k$.
	Since $\mathfrak{b}_\ell$ vanishes at $\eta$ and pulls back to $\mathfrak{I}^{\ell}$, $I_{\widetilde{\eta}}$ lies in the maximal ideal $\mathfrak{m}$ of this ring.
	By the threshold comparison in \cref{subsec:local-threshold-bounds}, monotonicity, and \cite[Lemma~5.7(2), (3)]{Singh-Takagi-Varbaro},
	\[
		\lct(Z;|L|) \leq \lct_\eta(Z;|L|) = \lct_{\widetilde{\eta}}(\widetilde{Z};I_{\widetilde{\eta}})
		\leq \lct_{\widetilde{\eta}}(\widetilde{Z};\mathfrak{m}) \leq k.
	\]

	For a cubic limit with $\theta > n-2$, each closed point $x \in X$ satisfies $\lct_x(X;|L|) \geq \lct(X;|L|) > n-2$, so $L$ is Cartier by \cref{cor:local-cartier-criterion}.
\end{proof}

\subsection{The full-dimensional image}
\label{subsec:full-dimensional-images}

For a full-dimensional image, its degree gives a lower bound for the last mixed intersection.
Combining this with \cref{lem:base-codimension-mixed-intersections,lem:codim_lct-intersection-bound} first bounds the codimension of the base locus, and then improves the threshold estimate.

\begin{proof}[Proof of \cref{prop:full-image-cartier}]
	If $\mathcal{B} = \varnothing$, then $L$ is Cartier since every non-Cartier point belongs to $\mathcal{B}$.
	Assume $\mathcal{B} \neq \varnothing$ and put $k = \codim_X \mathcal{B}$.

	The nondegenerate image $W\subset\PP^{n+1}$ is a hypersurface of degree at least two.
	If $d$ is the degree of $f:Y\to W$, then
	\[
		a_n=M^n=d\cdot\deg W\geq2.
	\]
	By \cref{lem:base-codimension-mixed-intersections}, every $a_i$ is at least two, while the first $k$ terms equal three.
	Therefore by \cref{lem:codim_lct-intersection-bound},
	\[
		k\geq\theta
		\geq\frac{n-1}{3(n+1)}
		\bigl(3k+2(n+1-k)\bigr),
	\]
	which yields
	\[
		k\geq\frac{n^2-1}{n+2}
		=n-2+\frac{3}{n+2}>n-2.
	\]
	Thus $k\geq n-1$.
	The first $n-1$ mixed intersections now equal three, and the last two are at least two.
	Substituting gives
	\[
		\theta\geq\frac{n-1}{3(n+1)}
		\bigl(3(n-1)+4\bigr)
		=n-2+\frac{n+5}{3(n+1)}>n-2.
	\]

	Hence $L$ is Cartier by \cref{lem:codim_lct-intersection-bound}.
\end{proof}

\section{The lower-dimensional image}
\label{sec:lower-image}
In this section, we show that the image $W$ of the rational map $\phi_{|L|}:X\dashrightarrow\PP^{n+1}$ must have dimension $n$.
\begin{proposition}
	\label{prop:full-dim-image}
	If $(X,L)$ is a K-semistable limit of smooth cubic hypersurfaces of dimension $n\geq 5$, then the image $W$ of $\phi_{|L|}$ has dimension $n$.
\end{proposition}
\begin{proof}
This is a consequence of \cref{prop:bounded-lower-image-reduction} and \cref{cor:full-dim-image-n-567}.
\end{proof}

Throughout this section, unless otherwise stated, let $Z$ be an $n$-dimensional K-semistable $\Q$-Fano variety with $n\geq5$, and let $L$ be an ample integral $\Q$-Cartier Weil divisor such that
\[
	-K_Z\sim_{\Q}(n-1)L,\qquad h^0(Z,L)=n+2.
\]
We retain the notation of \cref{subsec:product-ideal-model}; in particular, $W\subset\PP^{n+1}$ and $m=\dim W$.
By \cref{ext:higher-alpha},
\[
	\alpha(Z;L)\geq\frac{n-1}{n+1}>\frac12.
\]
Every cubic limit pair $(X,L)$ satisfies these assumptions by \cref{thm:cubic-limits}.

\subsection{The Stein factorization and adjunction}
\label{subsec:stein-factorization}

Let $p:Y\to Z$ and $f:Y\to W$ be the normalized graph of \cref{subsec:product-ideal-model}.
Factor $f$ as $Y\xrightarrow{h}B\xrightarrow{\beta}W$, where $h$ has connected fibers and $\beta$ is finite.
The Stein base $B$ is normal.
Put $A=\beta^*A_W$, an ample globally generated Cartier divisor.
Choose a smooth resolution $\tau:B'\to B$, and resolve the component of $Y\times_B B'$ dominating $B'$ to obtain a smooth $Y'$.
We have a commutative diagram
\[
	\begin{tikzcd}[row sep=large, column sep=large]
		Y' \arrow[d, "\sigma"'] \arrow[r, "h'"]
		& B' \arrow[d, "\tau"] & & \\
		Y \arrow[d, "p"'] \arrow[r, "h"]
		& B \arrow[r, "\beta"] & W \arrow[r, hook] & \PP^{n+1} \\
		Z \arrow[urr, dashed, bend right=12, "\phi"'] & & &
	\end{tikzcd}
\]
The morphism $h'$ and $\tau$ have connected fibers, as they are birational proper morphisms of normal varieties.
Set
\[
	A'=\tau^*A,\qquad p'=p\circ\sigma,\qquad
	f'=f\circ\sigma=\beta\circ\tau\circ h'.
\]
The divisor $A'$ is nef and big.
Using the pullback convention of \cref{subsec:product-ideal-model}, on $Y'$ we have
\[
	P=(p')^*L,\qquad M=(h')^*A'=(f')^*A_W,\qquad F=P-M.
\]
On the normalized graph, write $\overline{M}=p_*M$ for the movable Weil divisor class on $Z$.
The projection formula gives
\[
	\overline{M}=(p')_*M,\qquad
	L-\overline{M}=(p')_*F\geq0.
\]

Completeness of $|L|$ gives
\[
	H^0(Z,L)\cong H^0(Y,h^* A)\cong H^0(B,A)
	\cong H^0(B',A')\cong H^0(Y',M).
\]
The first identification is part of \cref{subsec:product-ideal-model}; the others follow from connected fibers and normality.
Thus $\beta:B\to W$ is the map defined by the complete system $|A|$.
We will use this to prove that $\beta$ has degree one.

These constructions use only the general hypotheses at the beginning of \cref{sec:full-image}.
The K-semistable Fano assumptions enter through the following semistability theorem.

\begin{theorem}[{\cite[Theorem~6(i) and Remark~4]{Druel-Guenancia-Paun}}]
	\label{ext:DGP-tangent}
	Let $Z$ be a complex K-semistable $\Q$-Fano variety.
	Then its reflexive tangent sheaf $T_Z$ is slope semistable with respect to $-K_Z$.
\end{theorem}

Since $-K_Z \sim_{\Q} (n-1)L$, this also gives slope semistability with respect to $L$.
We will show that $K_{B'}+(m-1)A'$ is not pseudoeffective.
The following adjunction theorem then classifies the polarized birational type of the Stein base.

\begin{theorem}[{\cite[Theorem~5.1(3) and Proposition~3.5(A)]{Andreatta}}]
	\label{ext:Andreatta}
	Let $(B,A)$ be a quasi-polarized variety of dimension $m\geq2$ with terminal $\Q$-factorial singularities, where $A$ is a nef and big Cartier divisor.
	Then $K_B+(m-1)A$ is not pseudoeffective if and only if any zero reduction of $(B,A)$ is one of $\bigl(\PP^m,\OK(1)\bigr)$, a quadric with its hyperplane polarization, a generalized cone over $\bigl(\PP^2,\OK(2)\bigr)$,
	or a $\PP^{m-1}$-bundle over a smooth curve with polarization restricting to $\OK(1)$ on the fibers.
	Here zero reductions and polarized birational equivalence are as in \cref{def:polarized-models}.
\end{theorem}

\subsection{The Stein base and its fibers}
\label{subsec:stein-base}

We first establish the adjunction obstruction.
Then we use the classification of \cref{ext:Andreatta} to show that the Stein base is a rational normal scroll.

\begin{lemma}
	\label{lem:adjunction-obstruction-for-the-stein-base}
	Let $Z$ be an $n$-dimensional K-semistable $\Q$-Fano variety with $n\geq2$, and let $L$ be an ample integral $\Q$-Cartier Weil divisor with $H^0(Z,L)\neq0$.
	Assume $-K_Z\sim_{\Q}(n-1)L$ and $0<m<n$.
	For the resolved Stein base $(B',A')$ of \cref{subsec:stein-factorization}, the divisor $K_{B'}+(m-1)A'$ is not pseudoeffective.
\end{lemma}

\begin{proof}
	If $K_{B'} + tA'$ is pseudoeffective for $t \geq 0$, then $(h')^*(K_{B'} + tA') \cdot P^{n-1} \geq 0$ since $P$ is nef.
	Using \cref{lem:kernel-foliation-bound} for rational map from $Z$ to $B'$ and $P-M = F \geq 0$, we obtain
	\[
		\frac{m(n-1)}{n} L^n
		\leq (h')^*(-K_{B'}) \cdot P^{n-1}
		\leq t M \cdot P^{n-1} \leq t L^n.
	\]
	Thus $t \geq m(n-1)/n > m-1$.
\end{proof}

\begin{proposition}
	\label{prop:classification-stein-base}
	Assume $m<n$.
	Let $\beta\colon B\to W$ be the finite morphism in the Stein factorization of a resolution of $\phi_{|L|}$.
	Then $m \geq 2$ and $\beta$ is an isomorphism.
	The projective image $W$ of $|L|$ is a rational normal scroll, which admits a canonical resolution 
	\[
		S=\PP_{\PP^1}\bigl(\OK(1)^{\oplus d}\oplus\OK^{\oplus(m-d)}\bigr),
		\qquad A_S=\OK_S(1), \qquad d = n-m+2 \leq m.
	\]
\end{proposition}

\begin{proof}
	Since $\alpha(Z;L)>1/2$, no member $D\in|L|$ contains a prime divisor $E$ on $Z$ with coefficient as least two.
	Otherwise $A_Z(E)=1$ would give $\alpha(Z;L)\leq\lct(Z;D)\leq1/2$.
	\begin{lemma}\label{lem:exclude-m-1}
		$\dim W =m \geq 2$.
	\end{lemma}
	\begin{proof}
		If $m=1$, then $h'\colon Y'\to B'$  is a surjective morphism to a smooth projective curve.
		Rationality of klt singularities identifies $H^1(Y',\OK_{Y'})$ with $H^1(Z,\OK_Z)$, and the latter vanishes by Kawamata--Viehweg vanishing.
		Hodge symmetry therefore gives $H^0(Y',\Omega_{Y'}^1)=0$.

		By the generic smoothness, the pullback map
		\[
			{h'}^*:H^0(B',\Omega_{B'}^1)\longrightarrow H^0(Y',\Omega_{Y'}^1)
		\]
		is injective.
		So $H^0(B',\Omega_{B'}^1)=0$ and $B'\simeq B \simeq \PP^1$ as $\tau$ has connected fibers.
		Then $h^0(B',A')=h^0(Z,L) = n+2$ gives $A'\cong \OK_{\PP^1}(n+1)$.
		For a general point $b$ of $B'$, its fiber of $h'$ intersects the open subset where $p'$ is an isomorphism.
		Since  $H^0(B',A'(-2b))\cong H^0(\PP^1,\OK_{\PP^1}(n-1))\neq 0$, pullback along the rational map $Z \dashrightarrow B'$, a section of $A'$ vanishes twice at that point yields a member of $|L|$ with a multiple prime component, contradicting the alpha bound.
		Also $m=0$ is impossible because the sections would all be proportional.

	\end{proof}

	For $m\geq2$, let $(B',A')$ be as in \cref{lem:adjunction-obstruction-for-the-stein-base}, we have a classification of poalrized birational models of $(B',A')$ by \cref{ext:Andreatta}.
	The projective space and quadric cases have respectively $m+1$ and $m+2$ independent sections, both smaller than $n+2$, and hence are excluded as polarized birational equivalence preserves every section space $H^0(kA')$ and the complete linear-series maps.

	The generalized Veronese-cone case is excluded by double hyperplanes arguments.
	A square of a linear form on $\PP^2$ gives a hyperplane section equal to twice the cone over its line.
	Vary the line and we get a moving family of hyperplane sections. 
	The pullback of the general hyperplane section of this family along the rational map gives a divisor in $|L|$ on $Z$ with coefficient at least two along some component.
	This again contradicts the alpha bound.

	It remains to consider the projective-bundle alternative.
	Its base curve is $\PP^1$ by the preceding argument in \cref{lem:exclude-m-1}.
	Up to tensoring a line bundle from the base, we may write the zero reduction of $(B',A') $ as $S=\PP_{\PP^1}\bigl(\bigoplus_{i=1}^m\OK(e_i)\bigr)$ with $A_S=\OK_S(1)$.
	The polarization is nef, so $e_i\geq0$.
	If some $e_i\geq2$, choose a section supported in that summand with a double zero at a general point of $\PP^1$.
	Its hyperplane divisor contains twice the fiber over that point.
	The same pullback argument contradicts $\alpha(Z;L)>1/2$.
	Therefore all $e_i$ belong to $\{0,1\}$.
	Let $d=\sum_{i=1}^m e_i$.
	The section count gives
	\[
		n+2=h^0(S,A_S)=m+d,\qquad d=n-m+2,\qquad d\leq m.
	\]

	It remains to identify the Stein base.
	Polarized birational equivalence identifies all section spaces of $B'$ and $S$, and therefore identifies their complete linear-series maps.
	Hence $W$ is the image of the birational tautological morphism  of $S$, which is a rational normal scroll of degree $d$ in $\PP^{n+1}$.
	As the complete map of $(B,A)$ also has degree one,  the morphism $\beta:B\to W$ is finite birational to a normal variety, and hence an isomorphism.
\end{proof}

Write $\varpi:S\to\PP^1$ for the bundle projection and put $R_{1,S} = \varpi^*\OK_{\PP^1}(1)$. 
We call the projective bundle $S$ the \emph{canonical resolution} of $W$ and denote the tautological morphism by $\rho:S\to W$.

The classification alone still permits scrolls of larger degree.
However, the two-ruling inequalities below provide a strong upper bound on the degree of the scroll.
This is the second use of tangent-sheaf semistability.

\begin{lemma}\label{lem:the-two-ruling-inequalities}
	Assume $m<n$.
	There are effective movable Weil divisor classes $R_1$ and $R_2$ on $Z$ such that $L-R_1-R_2$ is effective and
	\[
		\frac{R_1\cdot L^{n-1}}{L^n}\geq\frac{n-1}{2n},\qquad
		\frac{R_2\cdot L^{n-1}}{L^n}\geq\frac{(n-1)(d-1)}{nd}.
	\]
	In particular,
	\[
		d\leq 2+ \frac{4}{n-3},
	\]
	where $d=n-m+2$ is the degree of the scroll.
\end{lemma}

\begin{proof}
	The $d$ sections of $A_S-R_{1,S}$ from the positive summands define the second ruling $S \dashrightarrow \PP^{d-1}$.
	On a common resolution $\pi\colon \widehat{Y} \to Z$ mapping to $S$, let $q_i\colon \widehat{Y} \to \PP^{b_i}$ be the two induced morphisms, with $(b_1,b_2) = (1,d-1)$.
	Put $R_i = \pi_*q_i^*\OK_{\PP^{b_i}}(1)$. The two divisor classes on $S$ sum to $A_S$. Removing the fixed part on $\widehat{Y}$ and pushing forward gives $R_1 + R_2 \leq \overline{M} \leq L$.
	Applying \cref{lem:kernel-foliation-bound} to each $q_i$ gives the stated lower bounds, since $-K_{\PP^{b_i}} = (b_i+1)\OK_{\PP^{b_i}}(1)$.
	Their sum is at most one, so $\frac{n-1}{n}\bigl(\frac12 + \frac{d-1}{d}\bigr) \leq 1$, as claimed.
\end{proof}

\begin{proposition}
	\label{prop:bounded-lower-image-reduction}
	If $n \geq 8$, the complete linear system $|L|$ is generically finite onto its image.
	If $n \in \{5,6,7\}$ and $|L|$ is not generically finite, its image has dimension $n-1$ and degree $d=3$.
\end{proposition}

\begin{proof}
	Suppose $m < n$. By \cref{prop:classification-stein-base}, $d = n-m+2 \geq 3$ and $d \leq m$, so $2d \leq n+2$.
	If $n \geq 8$, \cref{lem:the-two-ruling-inequalities} gives $d < 3$, a contradiction.
	For $n = 5$, the section count gives $2d \leq 7$. For $n \in \{6,7\}$, the two-ruling bound gives $d < 4$.
	Thus $d = 3$ and $m = n-1$ in all three cases.
\end{proof}

The remaining images have dimension one less than $Z$.
In these cases, we next show that their general resolved fibers are rational curves.
The ruling of the rational normal scroll makes the induced rank-one foliation positive along a general complete-intersection curve, so we can apply the following theorem.

\begin{theorem}[{{\cite[Theorem~1]{Kebekus-SolaConde-Toma}}, see also \cite[Main theorem]{BM_Rational_2016}}]
	\label{ext:KSCT}
	Let $Z$ be a normal complex projective variety, let $\Gamma\subset Z_{\mathrm{reg}}$ be a complete curve, and let $\mathcal{F}\subset T_Z$ be a possibly singular foliation which is regular along $\Gamma$.
	If $\mathcal{F}|_{\Gamma}$ is an ample vector bundle, then the leaf through every point of $\Gamma$ is algebraic, and the closure of the leaf through a general point of $\Gamma$ is rationally connected.
\end{theorem}

Before we proceed, we recall the construction of \cref{subsec:stein-factorization} with $B'=S$.
\[
	\begin{tikzcd}[row sep=large, column sep=large]
		Y' \arrow[dr, "f'"] \arrow[d, "\sigma"'] \arrow[r, "h'"] \arrow[dd,  bend right=50, "p'"']
		& S \arrow[d, "\tau=\rho"] & & \\
		Y \arrow[d, "p"'] \arrow[r, "h=f"]
		&  W \arrow[r, hook] & \PP^{n+1} \\
		Z \arrow[ur, dashed, bend right=12, "\phi"'] & & &
	\end{tikzcd}
\]

\begin{lemma}\label{lem:residual-fibers-are-rational}
	If $m=n-1$, then a smooth connected general fiber of a resolution of $\phi_{|L|}$ is isomorphic to $\PP^1$.
\end{lemma}

\begin{proof}
	The projective-bundle formula gives $K_S = -(n-1)A_S + R_{1,S}$.
	By \cref{lem:kernel-foliation-bound}, the rank-one kernel foliation $\mathcal{F}$ on $Z$ induced by $h'$ satisfies
	\[
		c_1(\mathcal{F}) \cdot L^{n-1}
		\geq \bigl((n-1)(L-\overline{M}) + R_1\bigr) \cdot L^{n-1} > 0,
	\]
	as $L-\overline{M} \geq 0$ and $R_1 = (p')_*(h')^*R_{1,S}$ is nonzero and effective.

	Choose a sufficiently general complete-intersection curve $\Gamma \subset Z$ using a large very ample Cartier multiple of $L$.
	Since $\operatorname{Sing}(Z)$  and  $\operatorname{Sing}(\mathcal{F})$ both have codimension at least two, we can assume $\Gamma$ is smooth and avoids both singular loci.
	The restriction $\mathcal{F}|_{\Gamma}$ is a line bundle of positive degree, so it is ample.
	By \cref{ext:KSCT}, the closure of the leaf through a general point of $\Gamma$ is rationally connected.

	By generic smoothness, there is an irreducible open subset $U$ of $S$ such that the restriction of $h'$ is a smooth proper family of connected curves with constant genus over $U$. 
	Choose a nonempty open subset $V\subset (h')^{-1}(U)$ such that $p'|_V$ maps $V$ isomorphically onto an open subset $Z^\circ\subset Z_{\mathrm{reg}}$.
	Let $g=h'|_V\circ(p'|_V)^{-1}$ be the induced morphism $Z^\circ\to U$.
	The foliation $\mathcal{F}$ is tangent to the fibers of $g$.

	For a sufficiently general choice of $\Gamma$ above, take a general point $x \in \Gamma \cap Z^\circ$, and put
	$s = g(x) \in U$.
	The smooth connected fiber $C_s =(h')^{-1}(s) \subset Y'$ maps birationally under $p'$ onto the closure $\Lambda_x \subset Z$ of the leaf through $x$.
	Then \cref{ext:KSCT} shows that $\Lambda_x$ is rationally connected, so its normalization $C_s$ has genus zero.
	Constancy of the genus over $U$ now implies that every closed fiber of $(h')^{-1}(U) \to U$ is isomorphic to $\PP^1$.
	Since $\rho:S \to W$ is birational, the general fibers of $h'$ and $f' = \rho \circ h':Y' \to W$ coincide.
	Further smooth birational resolutions preserve this conclusion, since their smooth proper generic curve fibers are birational and hence isomorphic.
\end{proof}

\subsection{The remaining cases}
\label{subsec:residual-cases}

Throughout this subsection, we assume in addition that
\[
	L^n=3,\qquad 2L\text{ is Cartier},\qquad m<n.
\]
As in \cref{prop:bounded-lower-image-reduction}, we assume $n\in\{5,6,7\}$ and the Stein base is the rational normal scroll of degree $d=3$ and dimension $m=n-1$.

Taking $\ell = 2$ in \cref{subsec:product-ideal-model} shows that $2F$ is Cartier. Since $L^n = 3$ and $\lambda = n-1$, the coefficient in \cref{lem:mixed-intersection-threshold} is
\[
	\theta=\frac{n-1}{3(n+1)}\sum_{i=0}^n a_i.
\]

Let $C$ be a smooth connected general fiber of $f'$.
By \cref{lem:residual-fibers-are-rational}, one has $C\cong\PP^1$.
The degree of the general fiber $C$ with respect to $P$ has only two possibilities, characterized by the codimension of the base locus $\mathcal{B}$ of $L$.

\begin{lemma}\label{lem:fiber-degree-dichotomy}
	One has
	\[
		a=P\cdot C\in\Bigl\{\frac12,1\Bigr\}.
	\]
	Moreover, $\mathcal{B}$ is nonempty and, writing $k = \codim_Z\mathcal{B}$, one has
	\[
		a = \frac12 \Longleftrightarrow k \leq n-1,
		\qquad
		a = 1 \Longleftrightarrow k = n.
	\]
\end{lemma}

\begin{proof}
	Since $2P=(p')^*(2L)$ is Cartier, we have $2a\in\Z$.
	The birational morphism $p':Y'\to Z$ is an isomorphism on a dense open subset met by the general fiber $C$ of $f':Y'\to W$.
	Hence $(p')_*C$ is a nonzero effective curve, and ampleness of $L$ gives
	\[
		a=P\cdot C=L\cdot(p')_*C>0.
	\]
	Since $M\cdot C = 0$, one also has $F\cdot C = a > 0$, so $\mathcal{B}$ is nonempty.
	Set $k = \codim_Z\mathcal{B}$.

	Since $C$ is a smooth connected general fiber of $f'$, the projection formula yields
	\[
		a_{n-1}=PM^{n-1} =(f')^*A_W^{n-1} P =3(P\cdot C)=3a.
	\]
	If $k \leq n-1$, by \cref{lem:base-codimension-mixed-intersections} we have
	\[
		3 = a_0 = a_{k-1} > a_k \geq a_{n-1} = 3a.
	\]
	Thus $0 < a < 1$.
	Since we already know that $2a\in\Z$, it follows that $a = 1/2$.
	If $k = n$, the same lemma gives $a_{n-1} = a_0 = 3$, so $a = 1$.
	This completes the proof.
\end{proof}

It remains to exclude the case $P\cdot C=1/2$, or equivalently, to show that $\mathcal{B}$ is finite.
As in the full-dimensional case (proof of \cref{prop:full-image-cartier}), we first use the threshold estimate to obtain a lower bound for the codimension of $\mathcal{B}$, and then use this codimension bound to sharpen the threshold estimate.

\begin{lemma}
	\label{lem:half-degree-fiber-threshold}
	Suppose $P\cdot C = 1/2$.
	Then $\mathcal{B}$ has codimension at least $n-2$, and
	\[
		\theta \geq \frac{(n-1)^2}{n+1}
		 > n-3.
	\]
\end{lemma}

\begin{proof}
	Set $k = \codim_Z\mathcal{B}$. By \cref{lem:fiber-degree-dichotomy}, one has $k \leq n-1$ and
	\[
		a_{n-1} = PM^{n-1} = \frac32,\qquad a_n = 0.
	\]
	
	By \cref{lem:base-codimension-mixed-intersections}, one has $a_i \geq 3/2$ for $0 \leq i \leq n-1$ and $a_0 = \cdots = a_{k-1} = 3$.
	By \cref{lem:codim_lct-intersection-bound}, we also have
	\[
		k \geq \theta
		\geq \frac{n-1}{3(n+1)}
		\left(3k+\frac32(n-k)\right),
	\]
	whence
	\[
		k \geq \frac{n(n-1)}{n+3}
		= n-3+\frac{9-n}{n+3} > n-3.
	\]
	The last inequality uses $n \in \{5,6,7\}$.
	Thus $k \geq n-2$, and the first $n-2$ mixed intersections equal three.
	The two remaining terms through $a_{n-1}$ are at least $3/2$, so
	\[
		\theta \geq \frac{n-1}{3(n+1)}
		\left(3(n-2)+2\cdot\frac32\right)
		= \frac{(n-1)^2}{n+1}.
	\]
\end{proof}

The half-integrality of $F$ is decisive in the next argument: a general fiber $C$ with $P\cdot C= 1/2$ meets exactly one component $G$ of $F$, with $G\cdot C=1$.
We show that $G$ is also the only $p'$-exceptional divisor meeting $C$.
Adjunction along $C$ then gives a value for its log discrepancy that contradicts the preceding threshold estimate.

\begin{lemma}\label{lem:no-half-degree-fiber}
	The possibility $P\cdot C = 1/2$ does not occur.
\end{lemma}

\begin{proof}
	Suppose $P\cdot C = 1/2$.
	Take $C$ sufficiently general over the smooth locus of $W$, so that $f'$ is smooth near $C$.
	Write $F = \sum_j b_jD_j$,
	where the $D_j$ are distinct prime divisors and $b_j \in \tfrac12\Z_{>0}$, since $2F$ is Cartier.
	A general fiber is not contained in any $D_j$, so each $D_j\cdot C$ is a nonnegative integer.
	Since $M\cdot C = 0$, one has
	\[
		\frac12 = P\cdot C =F\cdot C = \sum_j b_j(D_j\cdot C).
	\]
	Consequently there is exactly one component $G$ of $\Supp F$ meeting $C$, and
	\[
		b_G = \frac12,\qquad G\cdot C = 1.
	\]
	By \cref{lem:half-degree-fiber-threshold}, the base locus has codimension at least $n-2 \geq 3$.
	Thus $|L|$ has no fixed divisorial component, and every component of $F$ is $p'$-exceptional.
	
	For a $p'$-exceptional prime divisor $H$ not in $\Supp F$, its center is not contained in the base locus.
	Near the generic point of that center, the map $\phi_{|L|}$ is a morphism, so $\dim f'(H) \leq \dim p'(H) \leq n-2$.
	Such a divisor is vertical over $W$ and does not meet $C$.
	
	Hence $G$ is the only $p'$-exceptional prime divisor meeting $C$.
	Write the discrepancy formula as
	\[
		K_{Y'}=(p')^*K_Z+\sum_{H\text{ exceptional}}(A_Z(H)-1)H.
	\]
	Since $f'$ is smooth near the general fiber $C$, its normal bundle is trivial and $K_{Y'}\cdot C=\deg K_C=-2$ as $C\simeq\PP^1$.
	Only $G$ contributes to the sum, and $G\cdot C=1$; hence
	\[
		A_Z(G)-1
		= (K_{Y'}-(p')^*K_Z)\cdot C
		= -2+\frac{n-1}{2}.
	\]
	Thus $A_Z(G) = (n-3)/2$.
	On the other hand, \cref{lem:mixed-intersection-threshold,lem:half-degree-fiber-threshold} give
	\[
		A_Z(G) \geq \theta b_G
		= \frac{\theta}{2} > \frac{n-3}{2},
	\]
	a contradiction.
\end{proof}

\begin{corollary}
	\label{cor:full-dim-image-n-567}
	Let $(X,L)$ be as in \cref{thm:cubic-limits} with $n \in \{5,6,7\}$. Then the image of $\phi_{|L|}$ has dimension $n$.
\end{corollary}

\begin{proof}
	Suppose $\dim W < n$. By \cref{prop:bounded-lower-image-reduction}, $W$ is a scroll of dimension $n-1$ and degree three.
	By \cref{lem:fiber-degree-dichotomy,lem:no-half-degree-fiber}, the base locus is nonempty and has codimension $n$.
	Thus \cref{lem:base-codimension-mixed-intersections} gives $a_0 = \cdots = a_{n-1} = 3$ and $a_n = 0$. Hence
	\[
		\theta = \frac{n(n-1)}{n+1} = n-2 + \frac{2}{n+1} > n-2.
	\]
	It follows from \cref{lem:codim_lct-intersection-bound} that $L$ is Cartier. It is then very ample by \cite[Corollary~1.7]{Fujita-del-Pezzo}, contradicting $\dim W < n$.
\end{proof}

\bibliographystyle{amsalpha}
\bibliography{kstability-cubics}

\end{document}